\def\editmode{0}

\def\bibfilenames{bibman_refs}

\def\spsformat{0}

\if\spsformat1

\documentclass{article}
\usepackage{spconf}
\usepackage{comment}

\else

\documentclass[conference,10pt]{IEEEtran}

\IEEEoverridecommandlockouts

\fi

\usepackage{savesym} 
\savesymbol{labelindent}
\usepackage{enumitem} 
\newcommand{\acom}[1]{\textcolor{red}{{[#1]}}} 

\if\editmode1  
\usepackage[backend=bibtex,style=alphabetic,sorting=debug,maxbibnames=99]{biblatex}
\DeclareFieldFormat{labelalpha}{\thefield{entrykey}}
\DeclareFieldFormat{extraalpha}{}
\bibliography{\bibfilenames}
\newcommand{\cmt}[1]{\noindent\textcolor{lightgreen}{\underline{[#1]}}} 
\newcommand{\hc}[1]{\textcolor{blue}{#1}} 

\setlistdepth{9}
\newlist{bulletlist}{enumerate}{9}
\setlist[bulletlist,1]{label=$\bullet$}
\setlist[bulletlist,2]{label=$\diamond$}
\setlist[bulletlist,3]{label=$\rightarrow$}
\setlist[bulletlist,4]{label=$\circ$}
\setlist[bulletlist,5]{label=$-$}
\setlist[bulletlist,6]{label=$\square$}
\setlist[bulletlist,7]{label=$\star$}
\setlist[bulletlist,8]{label=$\checkmark$}
\setlist[bulletlist,9]{label=$\Delta$}
\newenvironment{bullets}{\begin{bulletlist}}{\end{bulletlist}}

\newcommand{\blt}[1][noargpassed]{
  \item%
  \ifthenelse{\equal{#1}{noargpassed}}{}{\cmt{#1}}%
}

\else
\usepackage{cite}
\newcommand{\cmt}[1]{} 
\newcommand{\hc}[1]{\textcolor{black}{#1}} 
\newenvironment{bullets}{}{}
\newcommand{\blt}[1][noargpassed]{\ignorespaces}

\fi

\newcommand{\printmybibliography}{
\if\editmode1 
\printbibliography
\else
\bibliography{\bibfilenames}
\fi
}

\usepackage{fixltx2e}

\usepackage{graphicx}

\usepackage{subcaption}              

\usepackage[utf8]{inputenc}

\usepackage{amsfonts}
\usepackage{amsmath}
\usepackage{mathtools}

\usepackage{amssymb}

\usepackage{bm}

\usepackage{color,verbatim}
\usepackage{multirow}
\usepackage{accents}

\usepackage{theoremref}

\usepackage{url}

\newcounter{rulecounter}
\newcommand{\resetrule}{ \setcounter{rulecounter}{0}}
\resetrule

\newtheorem{myauxproblem}{Problem}

\newtheorem{myauxoptionalproblem}{Optional Problem}

\newsavebox{\selvestebox}
\newenvironment{colbox}[1]
  {\newcommand\colboxcolor{#1}%
   \begin{lrbox}{\selvestebox}%
   \begin{minipage}{\dimexpr\columnwidth-2\fboxsep\relax}}
  {\end{minipage}\end{lrbox}%
   \begin{center}
   \colorbox{\colboxcolor}{\usebox{\selvestebox}}
   \end{center}}

\definecolor{orange}{rgb}{1,0.8,0}
\definecolor{gray}{rgb}{.9,0.9,0.9}
\definecolor{darkgray}{rgb}{.3,0.3,0.3}
\definecolor{darkblue}{rgb}{.1,0.0,0.3}
\definecolor{lightblue}{rgb}{0.7,0.7,1}
\definecolor{lightred}{rgb}{1,0.7,.7}
\definecolor{purple}{RGB}{204,153,255}
\definecolor{lightgray}{rgb}{.95,0.95,0.95}
\definecolor{lightgreen}{rgb}{0.3,0.5,0.3}
\definecolor{darkgreen}{rgb}{0.05,0.3,0.05}

\newcommand{\ra}{$\rightarrow$~}

\newcommand{\tbm}[1]{{\tilde{\bm #1}}}

\newcommand{\cfield}{\mathbb{C}}
\newcommand{\rfield}{\mathbb{R}}

\newtheorem{myproposition}{Proposition}
\newtheorem{myremark}{Remark}
\newtheorem{myproblemstatement}{Problem Statement}
\newtheorem{mylemma}{Lemma}
\newtheorem{mytheorem}{Theorem}
\newtheorem{mydefinition}{Definition}
\newtheorem{mycorollary}{Corollary}

\usepackage{balance}
\usepackage{hyperref}
\usepackage[]{algorithm2e}
\RestyleAlgo{ruled}
\SetKwInput{KwInit}{Initialize}

\newcommand{\nextv}[1]{}

\begin{document}

\title{Implicit Channel Charting \\with Application to UAV-aided Localization
  \thanks{
    This work has been funded by the IKTPLUSS grant 311994 of the Research
    Council of Norway.}}

\if\spsformat1
\name{Pham Q. Viet and Daniel Romero\thanks{Thanks to XYZ agency for funding.}}
\address{Author Affiliation(s)}
\else
\author{
    \IEEEauthorblockN{
        Pham Q. Viet and Daniel Romero
    }\\
    Dept. of Information and Communication Technology,
    University of Agder, Grimstad, Norway.\\
    Email:\{viet.q.pham,daniel.romero\}@uia.no.
  }
\fi

\maketitle

\newcommand{\pomit}[1]{} 
\newcommand{\comit}[1]{\textcolor{black}{#1}} 
\newcommand{\todo}[1]{\textcolor{blue}{#1}} 

\newcommand{\va}{\bm{a}}
\newcommand{\vu}{\bm{u}}
\newcommand{\vx}{\bm{x}}
\newcommand{\vr}{\bm{r}}
\newcommand{\mA}{\bm{A}}
\newcommand{\vb}{\bm{b}}
\newcommand{\vg}{\bm{g}}
\newcommand{\vcsi}{\tbm{g}}
\newcommand{\vI}{\bm{I}}
\newcommand{\vf}{\bm{f}}
\newcommand{\vp}{\bm{p}}
\newcommand{\vtheta}{\bm{\theta}}
\newcommand{\vy}{\bm{y}}
\newcommand{\vw}{\bm{w}}
\newcommand{\vz}{\bm{z}}
\newcommand{\tir}{\tilde{r}}
\newcommand{\userind}{{\hc{i}}}
\newcommand{\userindaux}{{\hc{j}}}
\newcommand{\uavposind}{{\hc{n}}}
\newcommand{\changes}[1]{\textcolor{blue}{#1}}
\newcommand{\unsure}[1]{\textcolor{red}{#1}}

\begin{abstract}
  
  Traditional localization algorithms based on features such as time
  difference of arrival are impaired by non-line of sight propagation,
  which negatively affects the consistency that they expect among
  distance estimates. Instead, fingerprinting localization is robust
  to these propagation conditions but requires the costly collection of
  large data sets. To alleviate these limitations, the present paper
  capitalizes on the recently-proposed notion of channel charting to
  learn the geometry of the space that contains the channel state
  information (CSI) measurements collected by the nodes to be
  localized. The proposed algorithm utilizes a deep neural network
  that learns distances between pairs of nodes using their measured
  CSI. Unlike standard channel charting approaches, this algorithm
  directly works with the physical geometry and therefore only
  implicitly learns the geometry of the radio domain. Simulation
  results demonstrate that the proposed algorithm outperforms its
  competitors and allows accurate localization in emergency scenarios
  using an unmanned aerial vehicle.

\end{abstract}


\section{Introduction}
\label{sec:intro}
\begin{bullets}%
  \blt[localization in cellular networks]
  \begin{bullets}%
    \blt[its applications]Localization services play a central role
    in countless applications such as navigation, augmented reality,
    autonomous driving, wireless communications and emergency
    response to name a few. 
\nextv{    \begin{bullets}
      \blt[wireless]For example, location information can be used in
      mobile networks to improve accuracy in beam
      alignment and channel estimation \cite{xiao2022integrated}
      \blt[emergency response]or in natural and man-made disasters to
      monitor the environment or find
      survivors~\cite{savvides2001dynamic}.
    \end{bullets}}
    %
  \end{bullets}%
  \blt[existing localization techniques]%
  \begin{bullets}%
    \blt[pilots]Most localization systems rely on algorithms that
    provide location estimates based on pilot signals that are
    received from satellites or terrestrial transmitters.
    \blt[LOS\ra model-based]In case of line-of-sight (LOS) reception,
      \emph{model-based} approaches are typically pursued, where
      geometric principles are applied to estimate locations from
      distance and/or angle estimates obtained from channel features
      such as time of arrival, time difference of arrival, or angle of
      arrival. 
      \blt[NLOS\ra data driven \ra fingerprinting]In turn, when there
      is not LOS to a sufficient number of transmitters, as occurs
      indoors or in urban scenarios, data-driven approaches are
      preferred since the aforementioned distance or angle estimates
      become too inaccurate.
    \begin{bullets}%
      \blt[description]The most prominent example of this class of
      algorithms is \emph{fingerprinting}, which involves recording a
      set of channel state information (CSI) vectors measured at known
      locations; see~\cite{sobehy2020nearest} and references therein.
      \begin{bullets}
        \blt[k-nn]Location estimates can be obtained, for instance, by
        comparing the CSI observed by the node to be located with the
        entries of this data set and applying K-nearest neighbors.
        \blt[DNN]More sophisticated alternatives rely on deep
        neural networks (DNNs) to learn a mapping from
        CSI~\cite{arnold2018deep, arnold2019novel} or from
        preprocessed
        CSI~\cite{niitsoo2018convolutional,li2019massive,ferrand2020feature}
        into location estimates.
      \end{bullets}
      \blt[limitation]The main limitation of fingerprinting approaches
      stems from the need for large data sets, which are costly to
      acquire since each entry involves obtaining the position of a
      sensor either manually or by means of auxiliary localization
      systems, e.g. by using a robot.
    \end{bullets}
  \end{bullets}%

  \blt[ss channel charting]To alleviate the cost of data collection,
  \emph{channel charting}~\cite{studer2018charting} has been recently
  proposed.
  \begin{bullets}%
    %
    \blt[general principle of CC]The idea is to establish a connection
    between the geometry in the \emph{radio space} where (features of)
    the CSI vectors reside and the geographical geometry of the
    \emph{physical space} where the nodes to be located lie. The key
    assumption is that CSI vectors acquired at spatially near
    locations are similar to each other.
    %
    %
    \blt[approaches]      
    \begin{bullets}
      \blt[no DNN]Fig.~\ref{fig:ecc_workflow}a depicts the main steps
      in channel charting.  
      \begin{bullets}%
        \blt[description]%
        \begin{bullets}%
          \blt[Workflow]There, a dimensionality reduction algorithm
          assigns a point in 2D or 3D space to each input CSI vector
          in such a way that the distance between each pair of points
          is similar in some sense to the dissimilarity between the
          feature representations of the CSI vectors acquired at those
          points; see Sec.~\ref{sec:preliminary}. This mapping is
          referred to as a \emph{channel chart}.
          \blt[anchors \ra semi-supervised]The relative positions of
          the points it returns approximately correspond to the
          relative positions of the nodes in the physical space. If
          in addition there are enough anchor nodes, i.e. nodes whose positions
          are known,  \emph{semi-supervised}
          extensions~\cite{huang2019improving} can
          provide absolute position estimates.
        \end{bullets}
        %

        \blt[limitations] In early works on channel charting, 
        feature extraction and dissimilarity metrics are manually
        engineered by relying on physical principles and heuristic
        considerations. 
      \end{bullets}
      \blt[neural network-based channel charting]%
      \begin{bullets}
        \blt[description]To reduce the
      inaccuracies arising from these approaches, DNN-based
      alternatives learn one of these steps from data.         
        \begin{bullets}%
          \blt[learning dissimilarity, fixed features]For instance,
          \cite{studer2018charting} and \cite{ferrand2020triplet} fix
          the feature extraction step and learn the dissimilarity
          metric or correspondence between the CSI vectors and the
          channel chart.
          %
          %
          \blt[learning features, dissimilarity fixed]Conversely,
          \cite{bromley1993signature} learns the mapping from CSI to
          features while fixing the dissimilarity metric to be the
          Euclidean distance.
          %
        \end{bullets}%
        \blt[limitations]
        \begin{bullets}%
          \blt[arbitrary]In short, both approaches learn only part of
          the workflow. 
          \blt[explicit]Besides, the explicit construction of a
          channel chart is convenient in those applications where
          only relative positions are required, but bypassing such a
          step is naturally expected to result in improved
          localization performance when absolute positions are
          needed.           
        \end{bullets}
      \end{bullets}      
    \end{bullets}
  \end{bullets}%
  
  \blt[our work]
  \begin{bullets}%
    \blt[main contribution: implicit CC loc]Building upon these two
    observations, the present work proposes
    \begin{bullets}%
      \blt[description]\textit{implicit channel charting-based
        localization} (ICCL), where the radio geometry is learned from
      data without explicitly constructing a channel chart.
      \begin{bullets}%
        \blt[step 1: NN]In the first step, a DNN is used to predict
        the physical (or geographical) distances between
        nodes given the CSI that they measure.
        \blt[step 2: multilateration]In the second step, these
        distances are utilized in combination with the locations of
        anchor nodes to estimate the absolute positions of the nodes. 
      \end{bullets}
      \blt[relative to other schemes]
      \begin{bullets}%
        \blt[CC]Thus, unlike most channel charting schemes, ICCL is
        supervised and provides absolute location estimates.
        %
        \blt[model-based loc]Relative to model-based localization
        algorithms, the proposed scheme inherits the robustness of
        fingerprinting to non-LOS (NLOS) propagation.
        \blt[FP]%
        \begin{bullets}%
          \blt[standard FP \ra learn geom. ]As compared  to conventional
          fingerprinting, the proposed algorithm learns the radio
          geometry from data,
          \blt[neural FP]whereas relative to DNN-based fingerprinting,
          learning is heavily improved since the fact that distances
          are learned instead of absolute positions gives rise to a
          natural data augmentation effect, where the number of
          training examples is quadratic in the number of entries of
          the data set; cf. Sec.~\ref{sec:training}. 
          \blt[more measurements]Finally, the proposed scheme
          leverages CSI acquired by multiple nodes rather than from
          only one, which is expected to increase robustness to noise
          and reduce the size of the required data set.
        \end{bullets}
      \end{bullets}%
    \end{bullets}

    \blt[side contribution: loc. with UAV]
    \begin{bullets}
      \blt[no infrastructure]Although the proposed ICCL approach could
      be used with arbitrary forms of CSI, this paper focuses on a
      scenario where an unmanned aerial vehicle (UAV) is used to
      locate nodes on the ground. This is well motivated when no
      terrestrial infrastructure is operational because of a natural disaster or a
      military attack and when no global navigation satellite systems
      (GNSSs) can be used, e.g. because nodes lack the appropriate
      sensors or because the propagation environment precludes LOS
      propagation from the satellites.
      %
      %
      \blt[relation to other works]This constitutes another
      contribution of the paper since, to the best of our knowledge,
      \begin{bullets}%
        \blt[loc. with UAVs \ra robust to channel impairments](i)
        existing schemes for localization with UAVs rely on
        model-based algorithms and therefore are sensitive to NLOS
        conditions and other channel impairments, and 
        \blt[1st CC+UAVs](ii) no previous work has considered channel
        charting in setups involving UAVs.        
      \end{bullets}
      %
    \end{bullets}%
  \end{bullets}

    \blt[structure of the paper] This paper is organized as follows.
    \begin{bullets}%
        \blt[preliminary]After reviewing some relevant background in Sec.~\ref{sec:preliminary},
        \blt[system model] Sec.~\ref{sec:problem_formualation}
        formulates the problem.
        \blt[localization] ICCL is proposed next in Sec.~\ref{sec:ICCL}
        \blt[simulation results]and its performance is empirically
        assessed in Sec.~\ref{sec:simulation}.
        \blt[conclusion] Finally, Sec.~\ref{sec:conclusion} concludes
        the paper.
    \end{bullets}%

    \blt[notation]\emph{Notation.} 
    \begin{bullets}%
      \blt[scalar, vector, matrix] Lower and uppercase boldface
      letters denote column vectors and matrices, respectively.
      \blt[identity matrix]$\vI$ denotes the identity matrix of
      appropriate size.
      \blt[conjugate transpose]The conjugate
      transpose operator is $(.)^H$.
      \blt[complex Gaussian distribution] A circularly-symmetric
      complex Gaussian distribution with mean $\mu$ and variance
      $\sigma^2$ is represented as
      $\mathcal{CN}\left(\mu,\sigma^2\right)$.
      \blt[Frobenius norm]Finally, $\|.\|$ denotes the Euclidean norm.
    \end{bullets}%
\end{bullets}

\section{Channel Charting}
\label{sec:preliminary}
\begin{bullets}%
  \blt[overview]Channel charting was proposed in
  \cite{studer2018charting} as an \emph{unsupervised} alternative to
  algorithms such as fingerprinting, which suffer from high data
  acquisition costs. In this context, \emph{supervised} means that
  each entry of the data set is a pair of a CSI vector and the
  location at which it was acquired, whereas \emph{unsupervised} means
  that each entry of the data set contains just a CSI vector. The
  price to be paid is that plain channel charting just provides coarse
  information of the relative locations of the nodes. In some
  applications, this kind of information suffices to enhance network
  functionalities such as handover management, predictive radio
  resource allocation, and user tracking or
  pairing~\cite{ferrand2020triplet}.

  \blt[key idea]As indicated earlier, the core idea behind channel
  charting is that spatially close sensors are expected to measure
  similar CSI from the relevant transmitters. 
  %
  \blt[standard CC] To apply this principle, the key steps of channel
  charting are described next and summarized in
  Fig.~\ref{fig:ecc_workflow}a. Consider $M$ nodes located at positions
  $\{\vp_\userind\}_{\userind=1}^M\subset\mathbb{R}^{D}$, where $D$ equals $2$ or
  $3$.
  \begin{bullets}%
  \blt[steps]
  \begin{bullets}%
    \blt[feature extraction]First, 
    \begin{bullets}%
      \blt[input]the CSI vector $\vcsi_\userind\in\mathbb{C}^L$ acquired by
      the $\userind$-th node is mapped
      \blt[output] into a feature vector $\vf_\userind=\bm \phi(\vcsi_\userind)\in\mathbb{C}^{L^\prime}$.
      \blt[algorithms]For example, such a transformation may involve
      computing second-order moments, scaling, and transforming
      the result into the angular domain~\cite{studer2018charting}.
    \end{bullets}%
    \blt[dissimilarity metric]For each pair of nodes, say $(\userind,\userindaux)$, a
    dissimilarity metric
    $d_{\userind,\userindaux}=\delta(\bm \phi(\vcsi_\userind),\bm \phi(\vcsi_\userindaux))$ is
    subsequently computed. Ideally, function $\delta$ should be chosen
    so that its returned value resembles
    the physical distance between the locations of 
    these nodes as much as possible. However, this is not generally
    doable and, for example,~\cite{agostini2020channel} uses the
    so-called \emph{correlation matrix distance}
    whereas~\cite{studer2018charting} uses Euclidean distance.
    \blt[dimensionality reduction]In the next stage, a dimensionality
    reduction algorithm is applied to find $M$ points
    $\{\vz_\userind\}_{\userind=1}^M\subset\rfield^D$ in such a way that the
    distance between the $\userind$-th and the $\userindaux$-th point is ideally
    $d_{\userind,\userindaux}$ for all $\userind,\userindaux$. Sammon's mapping \cite{sammon1969mapping}
    can be used to this end, but other methods such as principal
    component analysis (PCA)
    \cite{pearson1901liii
    } and autoencoders have also been considered~\cite{studer2018charting}. 

  \end{bullets}%
  \blt[channel chart]The mapping from $\vcsi_\userind$ to $\vz_\userind$ constitutes
  the channel chart. The vectors $\vz_\userind$  are named
  \emph{pseudopositions} because they approximately preserve the
  \emph{relative} positions of the vectors $\vp_\userind$. For this reason,
  the quality  of a channel chart is typically quantified by ad-hoc
  metrics such as the \emph{trustworthiness} and \emph{continuity}
  \cite{venna2001neighborhood,kaski2003trustworthiness,vathy2013graph}.
  \nextv{\acom{expand?}}%
\end{bullets}%
\blt[semi-supervised]However, it is also possible to obtain absolute
location estimates with channel charting by relying on semi-supervised
learning~\cite{huang2019improving}.

\begin{figure}
  \centering
  \includegraphics[width=9cm]{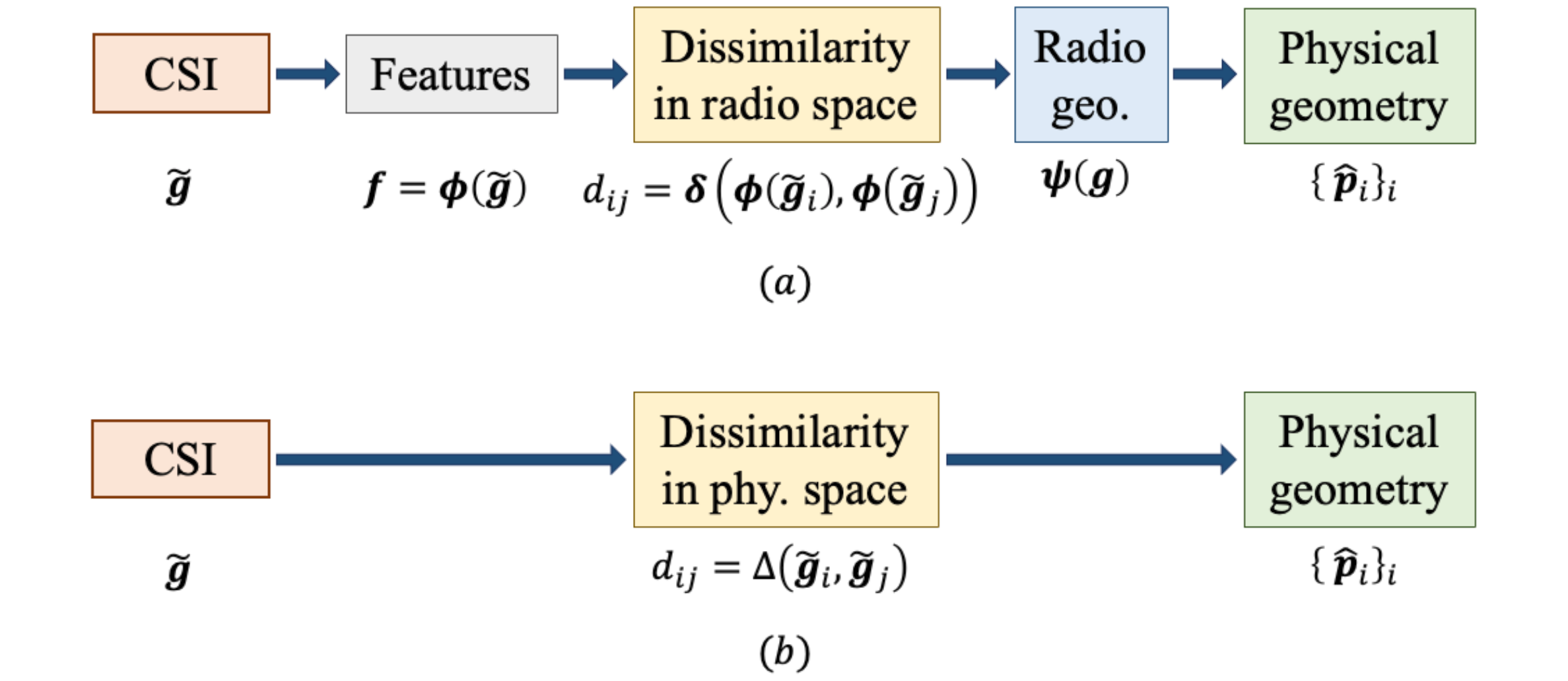}
  \caption{(a): Conventional (explicit) channel charting. (b):
    Implicity channel charting (proposed).}
  \label{fig:ecc_workflow}
\end{figure}    

%
\end{bullets}%

\section{Problem Formulation}
\label{sec:problem_formualation}

\begin{bullets}
  \blt[scenario description] Consider $M$ nodes located at positions
  $\{\vp_\userind\}_{\userind=1}^M\subset\mathbb{R}^{D}$, where $D$ equals $2$ or
  $3$. 
\begin{bullets}%
    \blt[nodes]
    \begin{bullets}%
      \blt[define anchors, unknowns]The positions
      $\mathcal{P}_{a}=\{\vp_{1},\vp_{2},\ldots,\vp_{M_a}\}$ of the
      first $M_a\geq 3$ nodes are known and, therefore, these nodes are
      referred to as \emph{anchors}. The locations
      $\mathcal{P}_{u}=\{\vp_{M_a+1},\vp_{M_a+2},\ldots,\vp_{M}\}$ of
      the rest of the nodes are unknown and, consequently, these nodes
      will be referred to as \emph{unknowns}.
      \blt[assumption of no localization technique]The unknowns are
      not able to localize themselves via GNSS, which occurs for
      example when (i) high buildings obstruct the LOS to satellites,
      (ii) the nodes are indoors, or (iii) the nodes are covered by
      debris, as occurs in applications where survivors from an
      earthquake must be located\nextv{\acom{cite?
          atif2021localization}}.  The unknowns cannot localize
      themselves using the terrestrial infrastructure either, which is
      relevant when the latter is not operational due to a natural
      disaster, a military attack, or a long
      blackout. \nextv{\acom{also distances between
          sensors?}}  
    \end{bullets}%
    
    \blt[uav]To localize the unknowns, a UAV flies over the area and
    transmits pilot signals at $N$ waypoints
    $\{\vu_{\uavposind}\}_{\uavposind=1}^N\subset\mathbb{R}^{3}$ along its trajectory.
    Although this paper considers a single UAV, it is straightforward
    to accommodate 
    multiple UAVs.
    \blt[CSI]For each of these $N$ waypoints, each node measures the
    CSI as described next. If the application at hand demands that the
    UAV locates the nodes, then all nodes report their measured CSI to
    the UAV. If, instead, each unknown must localize itself, the anchors
    send their measured CSI vectors to the UAV and the latter
    broadcasts them to all unknowns. 
    \blt[fig]The entire setup is illustrated in
    Fig.~\ref{fig:city_map}.

\begin{figure}
    \centering
    \includegraphics[width=7cm]{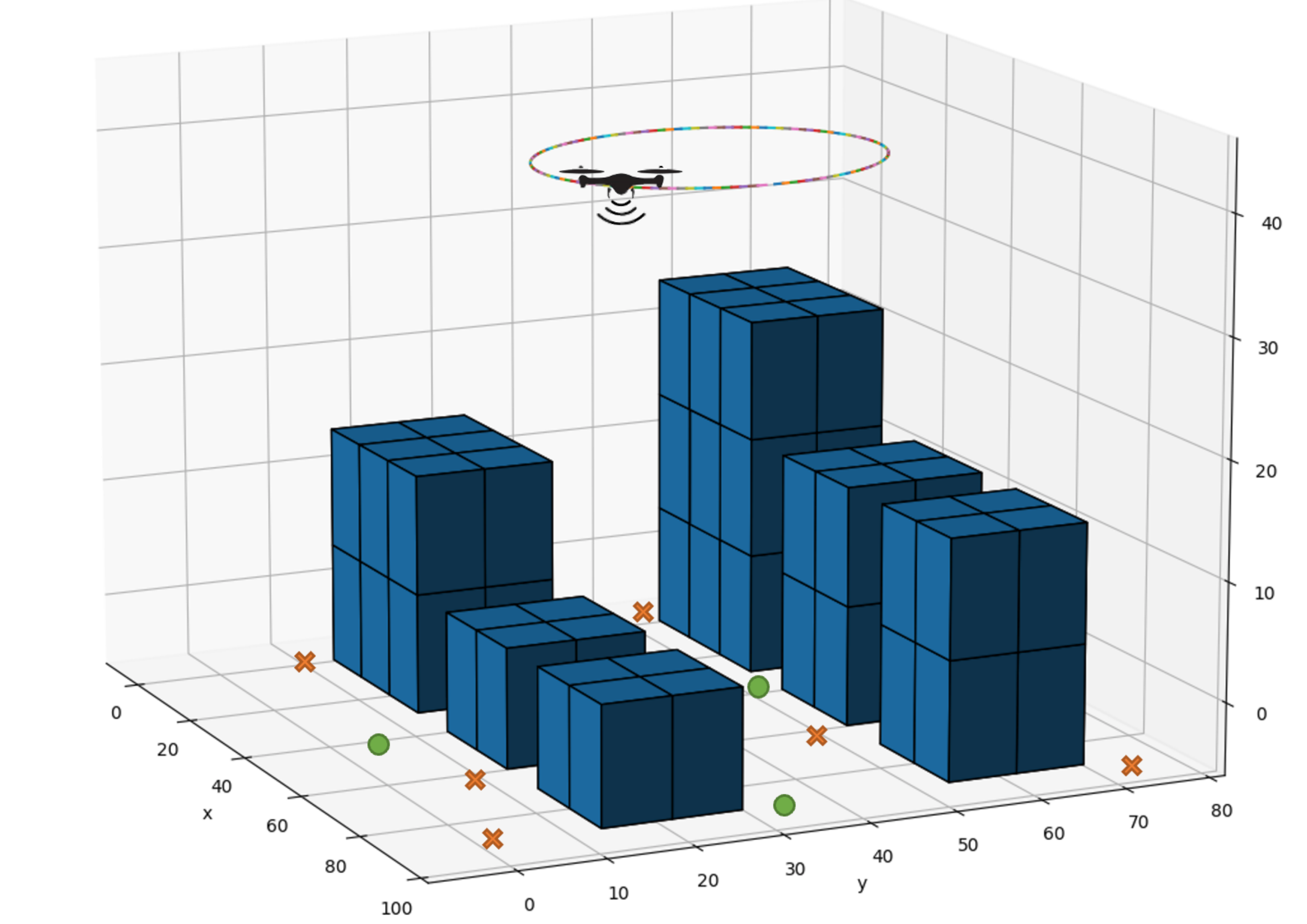}
    \caption{Illustration of a localization problem in an urban
      scenario using a UAV.  Green circles denote nodes with known
      locations. Orange crosses represent nodes with unknown
      locations. Blue blocks denote buildings.}
    \label{fig:city_map}
\end{figure}

\end{bullets}%
\blt[signal model]
\begin{bullets}%
  \blt[pilot sequence]At the $\uavposind$-th waypoint, the UAV transmits a
  pilot sequence consisting of $N_p$ symbols denoted as
  $ \vx_{\uavposind}=\left[x_{\uavposind}[1],x_{\uavposind}[2],\ldots,x_{\uavposind}[N_p]\right]^\top$.
    \blt[channel coefficients]For simplicity, assume that both the UAV
    and the nodes have a single antenna and that the channel is
    neither frequency nor time selective. Therefore, the channel
    between the $n$-th waypoint and the $\userind$-th node can be represented
    by a single coefficient $h_{\userind,\uavposind}\in \cfield$. The signal received
    at node $\userind$ is given by
    \begin{equation}
      \label{eq:rxsignal}
        \vy_{\userind,\uavposind} = h_{\userind,\uavposind}\vx_{\uavposind} + \vw_{\userind,\uavposind},
    \end{equation}
    where
    $\vw_{\userind,\uavposind}=\left[w_{\userind,\uavposind}[1],\ldots,w_{\userind,\uavposind}[N_p]\right]^\top$
    models  noise. 
\end{bullets}%

\blt[problem formulation] 
\begin{bullets}%
    \blt[given]Given
    \begin{bullets}%
        \blt[anchors]the anchor positions $\mathcal{P}_{a}$,
        \blt[transmitted + received signal]the pilot sequences
        $\{\vx_\uavposind\}_\uavposind$, and the received signals at all nodes
        $\{\vy_{\userind,\uavposind}\}_{\userind,\uavposind}$, 
    \end{bullets}%
    \blt[estimate]the problem is to estimate the  positions
    $\mathcal{P}_{u}$ of the unknowns.
\end{bullets}%
\end{bullets}

\section{Implicit Channel Charting-based Localization}
\label{sec:ICCL}
\begin{bullets}

  \blt[approach] This section proposes ICCL to solve the problem formulated in
  Sec.~\ref{sec:problem_formualation}. The algorithm consists of three
  phases.
\begin{bullets}
    \blt[CSI extraction] First, CSI needs to be extracted from the received signals.
    \blt[estimate distances] Given  the extracted CSI, a DNN
    predicts geographical distances between each pair of nodes.
    \blt[localization]Finally, the multilateration algorithm~\cite{savvides2001dynamic}  is used to recover
    the absolute positions of the unknowns given the aforementioned
    distances and the anchor locations.
\end{bullets}
The key steps in the proposed algorithm are shown in
Fig.~\ref{fig:ecc_workflow}b. Details of each phase will be provided in
the following subsections.
\end{bullets}

\subsection{CSI Extraction}
\label{sec:csiextraction}
\begin{bullets}%
  \blt[channel estimation]Although ICCL can be applied, in principle,
  to arbitrary forms of CSI, for concreteness and simplicity, CSI in
  this paper refers to the power gain. 

    \begin{bullets}%
        %
      \blt[channel estimate]In view of the model in \eqref{eq:rxsignal},
      the least-squares estimator of $h_{\userind,\uavposind}$ given $\vy_{\userind,\uavposind}$ and
      $\vx_{\uavposind}$ is given by
      \begin{equation}
        \label{eq:hestimator}
            \tilde{h}_{\userind,\uavposind} = {\vx_{\uavposind}^H\vy_{\userind,\uavposind}}/({\vx_{\uavposind}^H\vx_{\uavposind}}).
          \end{equation}
\nextv{         If $\vw_{\userind,\uavposind}\sim\mathcal{CN}\left(\bm 0,
           \sigma^2\vI\right)$, then  $\tilde{h}_{\userind,\uavposind}$ is also the
         \emph{minimum variance unbiased estimator}, the \emph{best
           linear unbiased estimator}, and the \emph{maximum
           likelihood estimator} of $h_{\userind,\uavposind}$~\cite{kay1}. In this case, observe
         that $\vy_{\userind,\uavposind}\sim\mathcal{CN}\left(\vx_{\uavposind}h_{\userind,\uavposind},
           \sigma^2\vI\right)$ and, as a result,}
          An estimate of the power gain can therefore be obtained as
          $\tilde{g}_{\userind,\uavposind}=|\tilde{h}_{\userind,\uavposind}|^2$. The CSI
          vector of the $\userind$-th node can then be defined as
          $\tilde{\vg}_\userind=\left[\tilde{g}_{\userind,1},\tilde{g}_{\userind,2},\ldots,\tilde{g}_{\userind,N}\right]^\top\in\mathbb{R}^{N}$.
          
          \blt[generation] For pre-training purposes, as discussed
          later, it is convenient to be able to generate samples of
          $\tilde{g}_{\userind,\uavposind}$ without simulating the propagation of the
          pilot signals through the channel as per
          \eqref{eq:rxsignal}. To this end, one can set
          $h_{\userind,\uavposind}=\sqrt{g_{\userind,\uavposind}}e^{j\varphi_{\userind,\uavposind}}$, where
          $g_{\userind,\uavposind}\in\rfield_+$ is the true power gain provided by
          some model, and
          $\varphi_{\userind,\uavposind}\sim\mathcal{U}\left(-\pi,\pi\right)$.
          Observe that if
          $\vw_{\userind,\uavposind}\sim\mathcal{CN}\left(\bm 0, \sigma^2\vI\right)$,
          then
          $\vy_{\userind,\uavposind}\sim\mathcal{CN}\left(\vx_{\uavposind}h_{\userind,\uavposind},
            \sigma^2\vI\right)$ and, as a result,
           $ \tilde{h}_{\userind,\uavposind}~\sim~\nextv{&\mathcal{CN}\left(\frac{\vx_{\uavposind}^H\vx_{\uavposind}
                h_{\userind,\uavposind}}{\vx_{\uavposind}^H\vx_{\uavposind}},\frac{\sigma^2}{\|\vx_{\uavposind}\|^2}\right)
            \\ &=}
          \mathcal{CN}\left(h_{\userind,\uavposind},{\sigma^2}/{\|\vx_{\uavposind}\|^2}\right)
          $.
          Then, one could equivalently write
          $\tilde{h}_{\userind,\uavposind}$ as
          $\tilde{h}_{\userind,\uavposind} =
          \sqrt{g_{\userind,\uavposind}}e^{j\varphi_{\userind,\uavposind}}
          + z_{\userind,\uavposind}$, where
          $z_{\userind,\uavposind}\sim\mathcal{CN}\left(0,{\sigma^2}/{\|\vx_{\uavposind}\|^2}\right)$
          models \emph{measurement error}. Since the noise is
          circularly symmetric, one can set
          $\varphi_{\userind,\uavposind}=0$ without loss of
          generality, which yields
          \begin{equation}
            \label{eq:ggen}
            \tilde{g}_{\userind,\uavposind}=\left|\tilde{h}_{\userind,\uavposind}\right|^2 = \left|\sqrt{g_{\userind,\uavposind}} + z_{\userind,\uavposind}\right|^2.
          \end{equation}
          Thus, samples of $\tilde{g}_{\userind,\uavposind}$ generated according to 
          \eqref{eq:ggen} are distributed as if the transmission of the
          pilot signals is simulated through \eqref{eq:rxsignal} and
          \eqref{eq:hestimator} is evaluated.
          
    \end{bullets}%
\end{bullets}%

\subsection{From CSI to Distances}
\begin{bullets}

  \blt[overview]This subsection presents the process of predicting
  distances between nodes from their CSI vectors $\{\tilde{\vg}_\userind\}_\userind$. A DNN is trained to
  this end and, therefore, it will be forced to implicitly learn the
  geometry in the CSI space.  

  \subsubsection{Architecture}
  \label{sec:architecture}
\begin{figure}
        \centering
        \includegraphics[width=7cm]{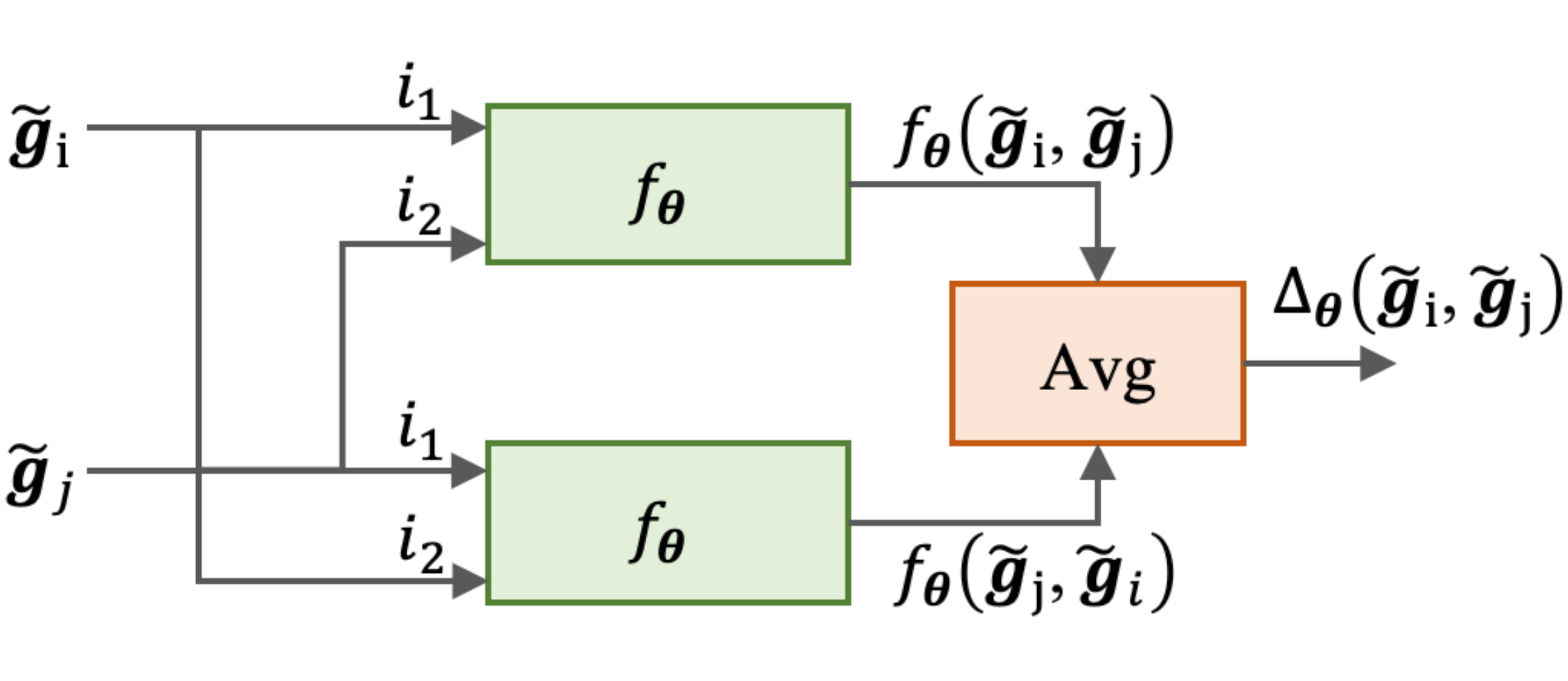}
    \caption{Architecture of the proposed network.\nextv{: (a) high level; (b) low level.}}
    \label{fig:architecture}
\end{figure}

\begin{bullets}
  \blt[input-output]Given the CSI vectors $\hat{\vg}_{\userind}$ and
  $\hat{\vg}_{\userindaux}$, the DNN obtains
  $\Delta_{\vtheta}(\tilde{\vg}_\userind,\tilde{\vg}_\userindaux)$, where $\vtheta$ is
  a vector collecting all its trainable parameters. This function will be
  fitted to the distances $\|\vp_\userind-\vp_\userindaux\|$. 
  \blt[symmetric high level] Since
  $\|\vp_\userind-\vp_\userindaux\| = \|\vp_\userindaux-\vp_\userind\| $, i.e., the distance from node
  $\userind$ to node $\userindaux$ equals the distance from node $\userindaux$ to node $\userind$, the
  learned function must be invariant to permutations of its inputs,
  i.e.
  $\Delta_{\vtheta}(\tilde{\vg}_\userind,\tilde{\vg}_\userindaux)=\Delta_{\vtheta}(\tilde{\vg}_\userindaux,\tilde{\vg}_\userind)$.
    \begin{bullets}%
      \blt[training] This could be approximately achieved while
      training by providing the network with each pair of nodes in
      both orders, i.e., with the examples
      $((\tilde{\vg}_\userind,\tilde{\vg}_\userindaux), \|\vp_\userind-\vp_\userindaux\|)$ and
      $((\tilde{\vg}_\userindaux,\tilde{\vg}_\userind), \|\vp_\userind-\vp_\userindaux\|)$ for all $\userind,\userindaux$.
      \blt[symmetric architecture]However, a more accurate and
      efficient approach is to impose invariance by means of the
      network architecture. To this end, one can let
        \begin{equation}
            \Delta_{\vtheta}(\tilde{\vg}_\userind,\tilde{\vg}_\userindaux) = \frac{1}{2}\left(f_{\vtheta}\left(\tilde{\vg}_\userind,\tilde{\vg}_\userindaux\right) + f_{\vtheta}\left(\tilde{\vg}_\userindaux,\tilde{\vg}_\userind\right)\right),
          \end{equation}
          where $f_{\vtheta}$ is a subnetwork. This is shown in
          Fig. \ref{fig:architecture}. Observe that, with this
          architecture, only the pairs of nodes with $\userind<\userindaux$ need to be
          provided at training time.
    \end{bullets}%
    
    \blt[detailed]\nextv{Fig. \ref{fig:low_level} provides the details
      of} For example, the subnetwork $f_{\vtheta}$ used in
    Sec.~\ref{sec:simulation} comprises the following layers:
    convolutional 2D, max pooling, convolutional 2D, max pooling,
    convolutional 2D, fully connected, and fully connected.
    \begin{bullets}%
      \blt[conv2D]Each 2D convolutional layer has 64 filters and
      $3\times 2$ kernels, except the last one, which
       has a $3\times 1$ kernel.
       \blt[maxpooling2d]The pool size of the 2D max-pooling layers is
       $2\times 1$.
       \blt[dense]The fully connected layers have 64 and 1 units,
       respectively.
    \end{bullets}%
\end{bullets}%

\end{bullets}%

\subsubsection{Training Process}
\label{sec:training}

\begin{bullets}%
  \blt[data collection]Training data can be collected in the same way
  as for fingerprinting. In the specific setup considered here, the
  UAV may start operating and sensors equipped with GNSS or other
  localization systems (e.g. as in LTE or 5G) can be sequentially
  placed at different positions where they measure the CSI. In case of
  emergency response applications, this measurement campaign is
  performed before the natural disaster or military attack.
    
  \blt[objective function]Once data is acquired, supervised learning
  is used to train the DNN. The cost function is the mean square error:
    \begin{equation*}
        C(\vtheta) \propto \nextv{=\frac{2}{M_0(M_0-1)}}\sum_{\userind=1}^{M_0-1}\sum_{\userindaux=\userind+1}^{M_0}\left[\Delta_{{\vtheta}}\left(\tilde{\vg}_\userind,\tilde{\vg}_\userindaux\right) - \|\vp_\userind-\vp_\userindaux\|\right]^{2},
      \end{equation*}
      where $M_0$ denotes the number of measurement locations in the
      data set.  Observe that the number of training examples is
      $M_0(M_0-1)/2$, whereas for DNN-based fingerprinting
      (cf. Sec.~\ref{sec:intro}) it would be just $M_0$. Thus, the DNN of
      ICCL is expected to be better trained than the DNN of DNN-based
      fingerprinting and, as a consequence, the former is expected to
      outperform the latter.
    
    %
    \blt[pretraining]Nonetheless, DNNs are known to be
    ``data-hungry''. Even with $M_0$ in the order of hundreds, $\bm
    \theta$ may not be learned properly if the network weights are
    initialized at random. Thus, it is convenient to pre-train the
    network using another data set, e.g. synthetically generated or
    measured in a different environment.

        
        
\end{bullets}%

\subsection{From Distances to Locations}

Given the distance estimates
$\hat d_{\userind,\userindaux} =
\Delta_{{\vtheta}}\left(\tilde{\vg}_\userind,\tilde{\vg}_\userindaux\right)$ provided
by the DNN as well as the anchor locations, ICCL estimates the
absolute positions of the unknowns via (possibly iterative)
multilateration \cite{savvides2001dynamic}. The possibility to use
this algorithm is a benefit of working directly with physical
distances rather than dissimilarity metrics in the radio geometry, as
in most channel charting algorithms.

\nextv{
\begin{bullets}
  \blt[algorithm]This paper uses multilateration
  \cite{savvides2001dynamic} to locate unknowns based on their
  distances to anchors.
    \blt[Given]Let $d_{i,j}\in\mathbb{R}_+$ be the predicted distance between anchor $i$ and unknown $j$. $\mathcal{D}=\{d_{i,j}\}$ is the set of all predicted distances between anchors and unknowns $i=1,\ldots,M_a;j=1,\ldots,M_u$.
    \blt[Estimating the $(n+1)$th location] 
    \begin{bullets}
        \blt[localizing] The following system of $n$ equations describes the geographical relationship between unknown $j$ and the anchors.
        \begin{equation}
            \left\{\begin{matrix}
            \|\vp_{j}-\vp_{1}\| = d_{j,1}\\
            \|\vp_{j}-\vp_{2}\| = d_{j,2}\\
            \vdots\\
            \|\vp_{j}-\vp_{M_a}\| = d_{j,M_a}
          \end{matrix}\right.\\
        \end{equation}
        After squaring and subtracting the last equation from the rest, the system of equation becomes the following least squares problem.
        \begin{equation}
            \label{eq:ls_problem}
            \min_{\vp_{j}} \|\mA_{j}\vp_{j} - \vb_{j}\|,
        \end{equation}
        where
        \begin{align}
          \label{eq:mA}
          \mA_{j} & = 2\begin{bmatrix}
            \vp_{M_a} - \vp_{1}, & \vp_{M_a} - \vp_{2}, & \ldots, & \vp_{M_a} - \vp_{M_a-1}
          \end{bmatrix}^\top\in\mathbb{R}^{(M_a-1)\times D};\\
          \label{eq:vb}
          \vb_{j} & = \begin{bmatrix}
            \|\vp_{M_a}\|^2 - \|\vp_{1}\|^2 - d_{M_a,j}^2 + d_{1,j}^2\\
            \|\vp_{M_a}\|^2 - \|\vp_{2}\|^2 - d_{M_a,j}^2 + d_{2,j}^2\\
            \vdots\\
            \|\vp_{M_a}\|^2 - \|\vp_{(M_a-1)}\|^2 - d_{M_a,j}^2 + d_{(M_a-1),j}^2
            \end{bmatrix}\in\mathbb{R}^{M_a-1}.
        \end{align}
        \blt[solution is position] Estimated location of unknown $j$ is the solution of LS problem \eqref{eq:ls_problem}, which is given by $\hat{\vp}_{j}=\left(\mA_{j}^\top\mA_{j}\right)^{-1}\mA_{j}^\top\vb_{j}$.
    \end{bullets}
\end{bullets}
}

\section{Experiments}
\label{sec:simulation}
\begin{bullets}%
    \blt[setup]
    \begin{bullets}%
      \blt[area] The simulation takes place in an urban area of size
      $100 \times 80$ m.
      \blt[CSI size] The UAV trajectory is a horizontal circle with
      center at (40, 45, 40)~m and radius 20~m.
      \blt[num CSI samples]At $N = 128$ waypoints, the UAV transmits a
      pilot signal with transmit power
      ${\|\vx_\uavposind\|^2}/{N_p}= 30~ \rm{dBm}$. However, the pilot signals
      are not explicitly generated; cf. Sec.~\ref{sec:csiextraction}. 
        \blt[data generation]
        \begin{bullets}%
            \blt[number + size of training samples]The CSI is then
            measured at $M_0=200$ positions drawn uniformly at random
            on the ground $(D=2)$.
            \blt[3D radio map]The true power gains $g_{\userind,\uavposind}$ are
            generated from the 3D city map depicted in
            Fig.~\ref{fig:city_map} using a tomographic
            model~\cite{patwari2008nesh} as in
            \cite{romero2022aerial}.
            \blt[true CSI]To focus on impairments in the testing
            phase, the noise power is set to 0 in the training data
            but it is greater than 0 for testing data. 
        \end{bullets}%
    \end{bullets}%
    \blt[compared algorithms] The proposed ICCL algorithm is compared
    with two algorithms. 
    \begin{bullets}%
      \blt[distance FP]One is the classical \emph{distance-based
        fingerprinting localization} (DFPL) algorithm;
      cf. Sec.~\ref{sec:intro}.
        \begin{bullets}%
          \blt[training data]This algorithm stores the training data.
          \blt[how it works]At testing time, given an input CSI
          vector, this method searches over the stored data and
          outputs the position that corresponds to the CSI vector that
          has lowest Euclidean distance to the input.
        \end{bullets}%
        \blt[neural FP]
        \begin{bullets}%
          \blt[network architecture]The second algorithm, termed
          \emph{neural-based fingerprinting localization} (NFPL), is
          similar in spirit to those in \cite{arnold2018deep,
            arnold2019novel,niitsoo2018convolutional,li2019massive,ferrand2020feature}
          but it is applied to the plain CSI vectors introduced in
          Sec.~\ref{sec:csiextraction}.  To obtain absolute position
          estimates, it trains a DNN with the same architecture as the
          subnetwork of ICCL (cf. Sec.~\ref{sec:architecture}) except
          for minor modifications to accommodate the different input
          and output size. Specifically, the kernels of the
          convolutional layers have size of $3\times 1$ instead of
          $3\times 2$ and the output layer has 2 neurons.
        \end{bullets}%
      \end{bullets}%
    \blt[pre-training]
    \begin{bullets}
      \blt[pretraining]Both NFPL and ICCL are pretrained with a data
      set that comprises $M_0=1000$ CSI vectors generated in a
      different environment, where the buildings have different
      dimensions.
        %
    \end{bullets}

    \blt[metric]To quantify the error between the true and estimated
    locations, the root mean square error 
\nextv{    \begin{equation}
      \textrm{RMSE}   = \sqrt{\frac{1}{M-M_a} \sum_{j=M_a+1}^{M}\mathbb{E}\left[\|\hat{\vp}_j-\vp_j\|^2\right]}
    \end{equation}
  }
  $ \textrm{RMSE} = [\frac{1}{M-M_a}
  \sum_{\userindaux=M_a+1}^{M}\mathbb{E}\left[\|\hat{\vp}_\userindaux-\vp_\userindaux\|^2\right]]^{1/2}$
  is used, where the expectation runs over realizations of the node
  locations and measurement noise.

    \blt[experiments]
    \begin{bullets}%
      \blt[rmse vs num. anchors]Fig. \ref{fig:rmse_anchors} shows the
      RMSE of ICCL vs. the number of anchors for different noise
      levels by
        \begin{figure}
            \centering
            \includegraphics[width=8cm]{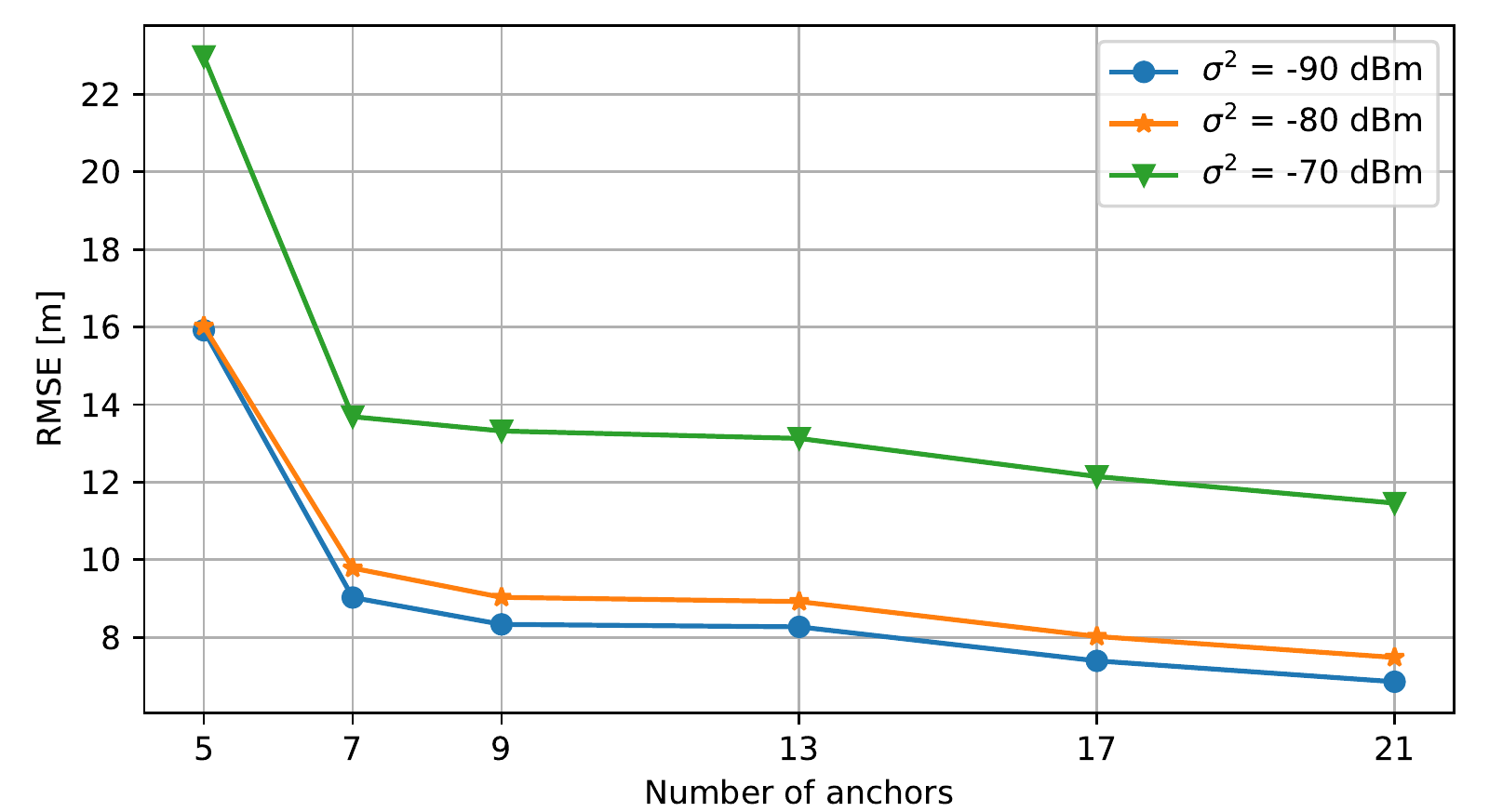}
            \caption{RMSE of the proposed ICCL algorithm.}
            \label{fig:rmse_anchors}
        \end{figure}
        \begin{bullets}%
          \blt[descriptions] averaging over 100 Monte Carlo realizations with
          $M=100$ nodes.
            \blt[comments]
            \begin{bullets}%
              \blt[more anchors, more precise] As expected, the
              more anchors, the more precise the estimated location.
              \blt[less than 10 meters] With only 7 anchors, the
              proposed algorithm can locate unknowns with less than
              10-meters average error provided that the noise power is
              sufficiently low.
            \end{bullets}%
          \end{bullets}%

          \blt[rmse vs noise]Fig.~\ref{fig:noise} depicts the RMSE of
          the compared algorithms vs. the noise level
        \begin{figure}
            \centering
            \includegraphics[width=8cm]{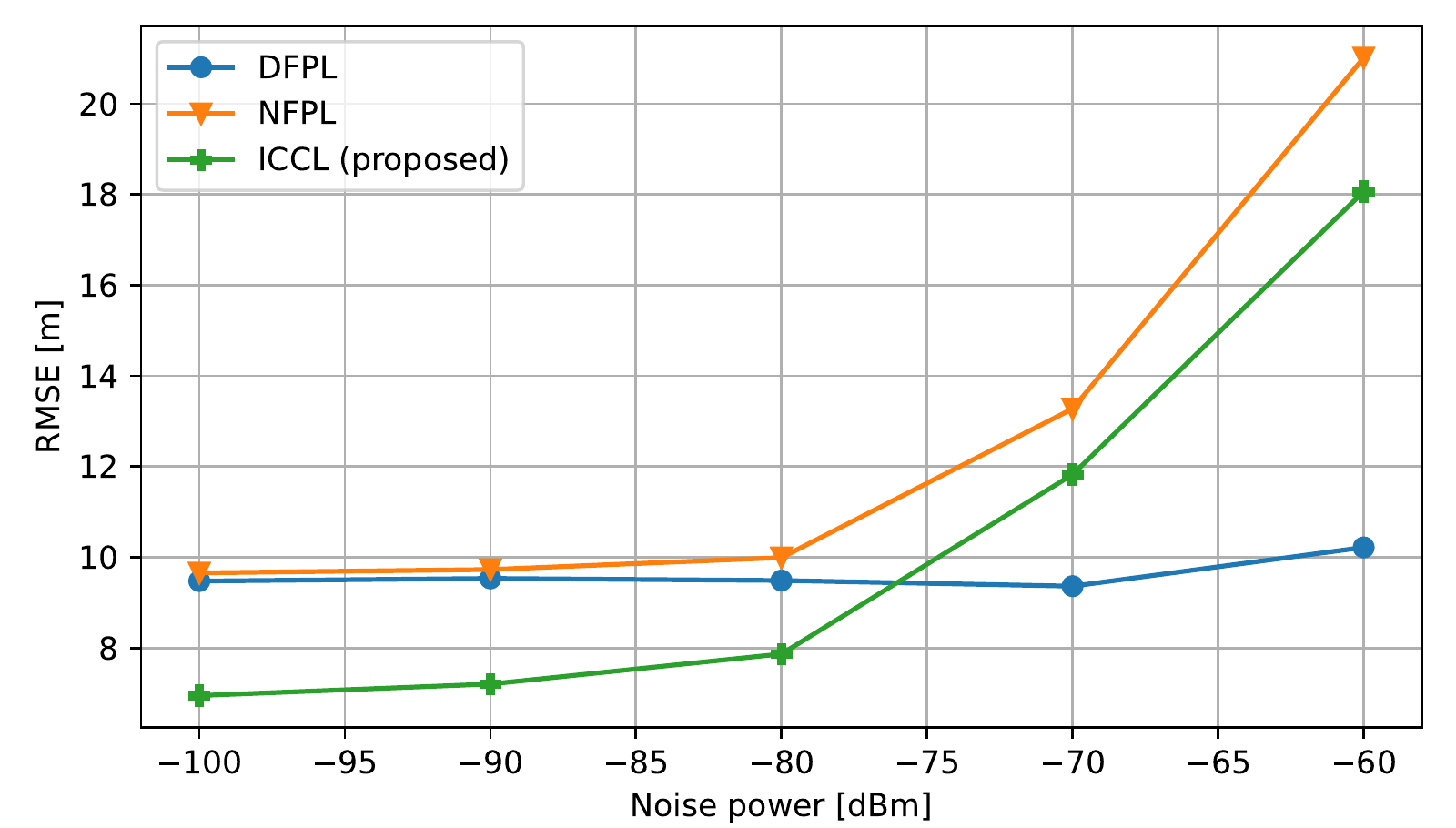}
            \caption{RMSE of the compared algorithms.}
            \label{fig:noise}
          \end{figure}
          \begin{bullets}%
            \blt[parameters]by averaging over  100 Monte Carlo
            realizations with $M_a= 20$ anchors and  $M-M_a=80$
            nodes. 
            \blt[comments]For a sufficiently small noise level, ICCL
            outperforms both DFPL and NFPL, which corroborates its
            ability to learn the radio geometry. However, at large
            noise power, the accuracy of the ICCL distance estimates
            degrades and DFPL works better. This is expected to
            improve if the training data is augmented by adding noise.
            \blt[FP]An apparently counterintuitive fact is that DFPL
            is seen to outperform NFPL. This phenomenon has already
            been observed in~\cite{sobehy2020nearest} and may be
            caused by the fact that DNNs require a large amount of
            training data. ICCL is less sensitive to this issue, as
            described in Sec.~\ref{sec:training}. In contrast, in
            \cite{ferrand2020feature}, NFPL offers a better
            performance than DFPL, but the reason may be that the
            latter applies a pre-processing step to the CSI
            vectors. Other works proposing NFPL schemes, such
            as~\cite{arnold2018deep,
              arnold2019novel,niitsoo2018convolutional,li2019massive},
            do not compare with DFPL.

          \end{bullets}%
        \end{bullets}%
\end{bullets}%

\section{Conclusions}
\label{sec:conclusion}
This paper proposes implicit channel charting-based localization
(ICCL) as a localization approach that implicitly learns the radio
geometry of a collection of CSI vectors from a data set. The idea is
inspired by channel charting and builds upon the well-known
fingerprinting localization method. Simulation results corroborate the
merits of the proposed approach.



\balance
\printmybibliography
\end{document}